\documentclass[12pt]{article}
\usepackage{amsfonts}
\usepackage{amsmath}
\usepackage{amssymb}
\def\be {\begin{equation}}
\def\ee {\end{equation}}
\def\ba {\begin{eqnarray}}
\def\ea {\end{eqnarray}}

\begin{document}
\title{ Distribution of Normalized Zero-Sets
of \\ Random Entire Functions with Small Random Perturbation}
\author{{ Weihong Yao}\\
{\small  Department of Mathematics, Shanghai Jiaotong University,}\\
{\small Shanghai 200240, P. R. China,}\\
{\small E-Mail: \ whyao@sjtu.edu.cn  }}
\date{}
\maketitle \baselineskip 18pt
\begin{center}
\begin{minipage}{120mm}
\textbf{A{\sc \textbf{bstract:}} } In this paper, we extend our
earlier result (see [Y-2008]) on the distribution of normalized zero-sets of random
entire functions to  random entire functions with small random perturbation.

\par
\
\end{minipage}
\end{center}

\section{ Introduction}
The well-known {\it fundamental theorem of algebra} states that
for every complex polynomial $P$, the degree of $P$ is equal to the
number of zeros of $P$ on the complex plane, counting multiplicities. This
suggests that one can use the counting functions (the number of zeros)
to measure the growth of $P$ (i.e. the degree of $P$). In 1929,
Nevanlinna extended the polynomial theory to meromorphic functions on ${\Bbb C}$
(which can be viewed as holomorphic maps $f: {\Bbb C}\rightarrow
{\Bbb P}^1({\Bbb C})$), in which case the growth function of $f$ is
given by its characteristic function $T_f(r)$ for
$|z|<r$. Geometrically, $T_f(r)$ is determined by the the area of the
image of $f(\bigtriangleup(r))$ in ${\mathbb P}^1({\Bbb C})$ with respect
to the Fubini-Study metric. Similar to the polynomials case,
Nevanlinna proved that  for almost all $a\in {\Bbb P}^1$, $N_f(r,
a)=T_f(r)+O(1)$ (or more precisely the integral average  of  $N_f(r, a)$
with respect to $a$ is $T_f(r)$). The result of this type is called the
{\it First Main Theorem of Nevanlinna}. Furthermore, Nevanlinna
obtained a much deeper result (called the  {\it Second Main Theorem of Nevanlinna})
which states that the sum of the difference $T_f(r)-N_f(r, a_j)$, for any
distinct $a_1, \dots, a_q\in  {\Bbb P}^1({\Bbb C})$, cannot exceed
$(2+\epsilon) T_f(r)$ for all $r\in (0, +\infty)$ except for a set of
finite measure. The theory is now known as {\it Nevanlinna theory} or
{\it value distribution theory}.  Nevanlinna's theory was later extended by H. Cartan and
L. Ahlfors to holomorphic curves.

\bigskip
The proof of the fundamental theorem of algebra comes from the
following observation: when we write $P(z)=a_nz^n+Q_{n-1}(z)$, where
$n=\deg P$, then $|Q_{n-1}(z)|<|a_nz^n|$ on $|z|=r$ for $r$ large
enough, hence Rouch\'e's theorem implies that the the zeros of $P$ is
the same as the zeros of $a_nz^n$. In other words, $P(z)$ can be
obtained from $a_nz^n$ through a small perturbation by $Q_{n-1}$.
Similarly, one can easily prove that the growth (characteristic
function) of $f$ is the same as $f+g$, the function obtained by
small perturbation by $g$. (Here, by small perturbations we mean
$T_g(r)=o(T_f(r))$).  Problems of these types are called {\it
  small perturbation problems} or called {\it problems of
  slowly moving targets}. In 1983, Steinmetz successfully extended
Nevanlinna's SMT to slowly moving targets, and in 1990, Ru-Stoll
extended H. Cartan's theorem to  slowly moving targets as well.

\bigskip
Recently, Shiffman and Zelditch, in their series of papers, initiated
the study of random value distribution theory.
The theory is based on
the following fundamental result of Hammersley:
{\it the zeros of random complex "Kac" polynomials
$f(z)=\sum_{j=0}^N a_jz^j$} ({\it where the coefficients $a_j$ are
independent complex Gaussian
random variables of mean 0 and variance 1}) tend to concentrate evenly about the unit circle
as the degree $N$ goes to the infinity.  Shiffman and Zelditch
extended the result to random polynomials of several complex
variables, as well as random holomorphic sections of line bundles. Their
method however largely relies on the use of kernel functions, thus
the strong ``orthonormal conditions'' are unavoidable.

\bigskip
In place of using kernel functions we propose a direct method
to study the value distribution of random meromorphic functions
(or maps).  This method starts with the computation of mathematical
expectations in the form of closed positive (1,1)-currents defined
by normalized
counting divisors, and is applicable to the much broader context
of random holomorphic functions and more generally
random meromorphic functions (maps).  As a first step in this
direction the author studied in [Y-2008] the
distribution of the normalized zero-sets of random holomorphic functions
$$G_n(z)=\sum_{\nu=0}^n\sum_{j_1=1}^\ell\cdots\sum_{j_\nu=1}^\ell
a_{j_1,\cdots,j_\nu}f_{j_1}(z)\cdots f_{j_\nu}(z),$$
where  $f_1(z),\cdots,f_\ell(z)$  are fixed holomorphic functions on a
domain $\Omega\in {\Bbb C}$, and coefficients $a_{j_1,\cdots,j_\nu}$ are independent complex
Gaussian random variables with mean $0$ and variance $1$.
More precisely, we studied the mean (mathematical expectation)
${\mathbf Z}\left(r,G_n\right)$
of the normalized counting divisor ${\mathbf Z}\left(r,G_n\right)$ of
$G_n(z)$ on the punctured disk
$0<\left|z\right|<r$ (in the sense of distribution) which is given
by
$${\mathbf Z}\left(r,G_n\right)=
\frac{1}{n}\sum_{G_n(z)=0,\atop
0<|z|<r}\left(\log\frac{r}{\left|z\right|}\right)\delta_z,
$$
where $\delta_z$ is the Dirac function.

\bigskip
We obtained the following result in [Y-2008].

\bigskip\noindent {\bf  Theorem A.} \ \ {\it Let $C$ be the smooth} ({\it possibly non-closed})
{\it curve in ${\mathbb C}$ consisting of all the points $z$ such that
$\left|f(z)\right|=\left(\sum_{j=1}^\ell\left|f_j(z)\right|^2\right)^{\frac{1}{2}}=1$
and
$f^\prime(z)=\left(f_1^\prime(z),...,f_l^\prime(z)\right)\not=0$.
Then the limit of  $ {\mathbf E}\left({\mathbf
Z}\left(r,G_n\right)\right)$ is equal to, in the sense of distribution,
$\log\frac{r}{\left|z\right|}$ times the sum of
$$\left|f(z)\right|\Xi\left(\left|f(z)\right|\right)\frac{\sqrt{-1}}{2\pi}\partial\bar\partial\log\left|f(z)\right|$$
and the measure on $C$ defined by the $1$-form
$$
\frac{\sqrt{-1}}{2}\sum_{j=1}^\ell\left(f_j(z)
d\overline{f_j(z)}-\overline{f_j(z)} df_j(z)\right),
$$
where
$$
\Xi\left(x\right)=\left\{\begin{matrix}\frac{2}{x}\ \ {\rm when\ \
}x>1\cr 1\ \ {\rm when\ \ }x=1\hfill\cr 0\ \ {\rm when\ \
}x<1.\hfill\cr\end{matrix}\right.
$$}

When $\ell =1$ and $f_1(z)=z$, our theorem recovers the result of
Hammersley.

\bigskip Theorem A can also be interpreted as an analogue of
Nevanlinna's First Main Theorem, which states that the integral average of the counting
function is equal to the growth (characteristic) function of $f$. We
also note that the approach used in [Y-2008] is very different from the method of
Shiffman-Zelditch.  In place of the use of sophisticated results
on kernel functions, we carry out a direct computation, which
basically comes down to the computation of the following
limit (see ``Complex Version of Lemma on the Convergence of integrals
as Distributions''), for $w\in{\mathbb C}^\ell$,
$$\lim_{n\to\infty}\frac{1}{n}\left(\frac{\sqrt{-1}}
{2\pi}\partial\bar\partial\log\left(\sum_{j=0}^n\left|w\right|^{2j}\right)\right),$$
in the sense of distribution.
In this paper, we extend Theorem $A$ to the random entire functions with
small perturbation, similar to the moving target case in Nevanlinna's
theory by an adaptation of the method and techniques of [Y-2008].
The main result of this paper is as follows: \\

\bigskip\noindent
{\it Main Theorem.}
Let $\Omega$ be an open subset of $\mathbb C$, $\ell$ be a positive integer
and $f_1(z),\ldots, f_\ell(z)$ be holomorphic functions on $\Omega$. For any $z_0\in\Omega$, let
$$
\kappa_{z_0,n},\lambda_{z_0,n},\xi_{z_0,n},\eta_{z_0,n}\ ,
$$
be four sequences of non-negative integers, $0\le n < \infty$. Assume that
$$
\lim\limits_{n\to\infty} \frac{\kappa_{z_0,n}}n = 0\ ,
\lim\limits_{n\to\infty} \frac{\lambda_{z_0,n}}n = 0\ ,
\lim\limits_{n\to\infty} \frac{\xi_{z_0,n}}n = 0\ ,
\lim\limits_{n\to\infty} \frac{\eta_{z_0,n}}n = 0\ .
$$
Let $A_z$ and $B_z$ be positive functions on $\Omega$.
For any nonnegative integer
$j$, let $g_j(z)$ be a holomorphic function on $\Omega$ such that

\begin{enumerate}
\item[(i)] $g_0(z)$ is nowhere zero on $\Omega$,

\item[(ii)] for each $0\le j<\infty$ and each point $j\in\Omega$
$$
\big|g_j(z)\big| \le (B_z)^{\kappa_{z,j}} \big(1+\big|f(z)\big|\big)^{\lambda_{z,j}}\ ,
$$
\item[(iii)] for each point $z\in\Omega$
$$
\liminf_{j\to\infty} (A_z)^{\xi_{z,j}} \big(1+\big|f(z)\big|\big)^{\eta_{z,j}} \big|g_j(z)\big|
$$
is positive.
\end{enumerate}

\noindent
Suppose furthermore that for each compact subset $K\subset \Omega$ there exists a positive
constant $C_K$ such that for each $z\in K$ and for each nonnegative integer
$j$ we have
$$B_z \le C_K; \quad \kappa_{z,j}, \lambda_{z,j}\le C_K $$

\noindent
For any positive integer $n$ let
$$
G_n(z) = \sum_{j=1}^n \sum_{1\le\nu_1,\ldots,\nu_j\le\ell} a_{\nu_1,\ldots,\nu_j}
g_j(z) f_{\nu_1}(z) \cdots f_{\nu_j}(z)
$$
be a random polynomial, where each coefficient $a_{\nu_1,\ldots,\nu_j}$ for
$1\le \nu_1,\ldots,\nu_j\le\ell$ and $0\le j\le n$ is an indeterminate which satisfies
the Gaussian distribution
$$
\frac 1\pi e^{-|z|^2}
$$
on $\mathbb C$, with the convention that $a_0$ is the single indeterminate
$a_{\nu_1,\ldots,\nu_j}$, when $j=0$. Let $\bold Z(G_n)$ be the {\it normalized counting divisor\/}
of $G_n(z)$ on $\Omega$ (in the sense of distribution) given by
$$
\bold Z(G_n) = \frac{\sqrt{-1}}n \sum_{G_n(z) = 0,\atop z\in\Omega} \delta_z\ ,
$$
where $\delta_z$ is the Dirac delta on $\mathbb C$ at the point $z$ of $\mathbb C$. Let
$\bold E(\bold Z(G_n))$ be the {\it expectation\/} of $\bold Z(G_n)$ which is defined by
$$
\int_{(a_{\nu_1,\ldots,\nu_j})\in\mathbb C^{N_{\ell,n}}}
\Big(\frac 1n \sum_{G_n(z)=0,\atop z\in\Omega} \delta_z\Big)
\prod\limits_{(a_{\nu_1,\ldots,\nu_j})\in\mathbb C^{N_{\ell,n}}}
\Big(\frac 1\pi\ e^{-|a_{\nu_1\cdots\nu_j}|^2} \frac{\sqrt{-1}}2\
da_{\nu_1,\ldots,\nu_j} \land \overline{da_{\nu_1\cdots\nu_j}}\Big)\ ,
$$
where $N_{\ell,n} = \frac{\ell^{n+1}-1}{\ell-1}$ is the number of indeterminates
$a_{\nu_1,\ldots,\nu_j}$ for $0\le j\le n$ and $1\le\nu_1,\ldots,\nu_j\le\ell$.
Let $C$ be the smooth (possibly non-closed) curve in $\Omega$ consisting of all the points
$z$ of $\Omega$ such that $\big|f(z)\big| = \Big(\sum_{j=1}^\ell |f_j(z)|^2\Big)^{\frac 12}=1$
and $f'(z)=(f'_1(z),\ldots,f'_1(z)) \ne 0$. Then $\bold E(\bold Z(G_n))$ is equal to the sum of
$$
\big|f(z)\big| \Xi \big(|f(z)|\big) \frac{\sqrt{-1}}{2\pi} \partial\bar\partial
\log \big|f(z)\big|
$$
and the measure on $C$ defined by the 1-form
$$
\frac{\sqrt{-1}}2 \sum_{j=1}^\ell \Big(f_j(z) \overline{df_j(z)} - \overline{f_j(z)} df_j(z)\Big)\ ,
$$
where
$$
\Xi(x) =
\begin{cases}
\frac 2x & \text{when $x>1$}\\
1        & \text{when $x=1$}\\
0        & \text{when $x<1$}\ .
\end{cases}
$$

\medskip\noindent
{\bf Remarks}
In the formulation of the Main Theorem, by (iii) we impose at each point
$z \in \Omega$ an asymptotic pointwise condition on a lower bound on $|g_j(z)|$.
It is natural to impose a condition in terms of lim inf, given that we want to allow
the holomorphic functions $g_j$ to have zeros.  On the other hand, in the computation
of mathematical expectations of normalized counting divisors, some
uniformity on compact subsets is required on asymptotic upper bounds on $|g_j(z)|$
in order to prove convergence of positive (1,1)-currents.

\medskip
\section{Proof of the Main Theorem}
We first recall the following key lemmas in [Y-2008].
\bigbreak\noindent(2.1)  {\bf Proposition [Y-2008] (Complex
Version of Lemma on the Convergence of Integrals as Distributions).}
{\it $$
\displaylines{\lim_{n\to\infty}\frac{1}{n}\left(\frac{\sqrt{-1}}{2\pi}\partial\bar\partial\log\left(\sum_{j=0}^n\left|z\right|^{2j}
\right)\right)\cr
=r\Xi\left(r\right)\frac{\sqrt{-1}}{2\pi}\partial\bar\partial\log r+
\left[\delta_{{}_{S^{2\ell-1}}}\right]\wedge
\frac{1}{r^2}\frac{\sqrt{-1}}{2}\sum_{j=1}^\ell\left(z_j d\bar
z_j-\bar z_j dz_j\right),\cr}
$$
is true, where $\left[\delta_{{}_{S^{2\ell-1}}}\right]$ denotes the
$1$-current on ${\mathbb C}^\ell$ defined by integration over
$$S^{2\ell-1}=\left\{\,z\in{\mathbb C^\ell}\,\left|\,|z|=\left(\sum_{j=1}^\ell\left|z_j\right|^2\right)^{\frac{1}{2}}=1\,\right.\right\}.$$}

By pulling back through
$f(z)=(f_1(z),f_2(z),\cdots,f_\ell(z)):\Omega\to{\mathbb C^\ell}$, the
above proposition implies

\bigskip
\noindent(2.2) {\bf Proposition [Y-2008].} \ \ {\it Let
$\Omega$ be a connected open subset of ${\mathbb C}$ and
$f(z)=(f_1(z),f_2(z),\cdots,f_\ell(z)):\Omega\to{\mathbb C^\ell}$ be
a nonconstant holomorphic function on $\Omega$. Let $C$ be the
smooth (possibly non-closed) curve in $\Omega$ consisting of all the
points $z$ of $\Omega$ such that $\left|f(z)\right|=1$ and
$f^\prime(z)=(f^\prime_1(z),f^\prime_2(z),\cdots,f^\prime_\ell(z))\not=0$,
where
$\left|f(z)\right|=(|f_1(z)|^2+|f_2(z)|^2+\cdots+|f_\ell(z)|^2)^\frac{1}{2}$.
Then the following
$$
\displaylines{\lim_{n\to\infty}\frac{1}{n}\left(\frac{\sqrt{-1}}{2\pi}\partial\bar\partial\log\left(\sum_{j=0}^n\left|f(z)\right|^{2j}
\right)\right)
=\left|f(z)\right|\Xi\left(\left|f(z)\right|\right)\frac{\sqrt{-1}}{2\pi}\partial\bar\partial\log
\left|f(z)\right|\cr + \left[\delta_{{}_{S^{2\ell-1}}}\right]\wedge
\frac{1}{\left|f(z)\right|^2}\frac{\sqrt{-1}}{2}\sum_{j=1}^\ell\left(f_j(z)
d\bar f_j(z)-\bar f_j(z) df_j(z)\right),\cr}
$$
is true, where $\left[\delta_{{}_{S^{2\ell-1}}}\right]$ denotes the
$1$-current on ${\mathbb C}^\ell$ defined by integration over
$$S^{2\ell-1}=\left\{\,f(z)\in{\mathbb C^\ell}\,\left|\,|f(z)|=\left(\sum_{j=1}^\ell\left|f_j(z)\right|^2\right)^{\frac{1}{2}}=1\,\right.\right\}.$$}

\bigskip
\noindent(2.3)
For the proof of the Main Theorem we will need to formulate a lemma on limits
of certain potential functions. To start with define on the domain $\Omega \subset \mathbb C$
the following subharmonic functions
$$\gamma_n = \dfrac 1{\, n\,} \log (1 + |f|^2 + \cdots + |f|^{2n}).$$
Define $\varphi : \Omega\to{\mathbb R}$ by
$$\begin{cases}
\varphi (z) = \log |f(z)|^2 & \mbox{ if } \ |f(z)|\ge 1\\
\varphi (z) = 0  & \mbox{ if }\  |f(z)|\le 1
\end{cases}
\, .
$$
In other words, $\varphi(z) = {\rm max}(0,\log|f|^2) = \log^+|f|^2$.
Then, we have

\smallskip
\noindent
{\bf Lemma 1}\
{\it
$\gamma_n (z)$ converges uniformly to $\varphi (z)$ on $\Omega$.
As a consequence,
$\sqrt{-1} \partial \overline \partial \gamma_n$ converges to
$\sqrt{-1} \partial \overline \partial \varphi$ as positive $(1,1)$-currents on
$\Omega$.
}

\smallskip
\noindent
{\it Proof.}\
For each positive integer $n$ define $\lambda_n: [0 , \infty) \to
{\mathbb R}$ by
$$\lambda_n (t)= \dfrac 1{\, n\,} \log (1+t+ \cdots + t^n).$$
For $0 \le t \le 1$ we have
$$ 0 \le \lambda_n (t) \le \dfrac 1{\, n\,} \log (n+1).$$
On the other hand, for $t\ge 1$ we have
$$ \log t = \dfrac 1{\, n\,} \log(t^n) \le \lambda_n (t) \le \dfrac 1{\, n\,} \log  \big((n+1)t^n\big ) =
\dfrac 1{\, n\,} \log (n+1) + \log t.$$
Let $\lambda : [0, \infty) \to {\mathbb R}$ be the monotonically
increasing continuous function defined by
$$
\begin{cases}
\lambda (t) = \log t & \mbox{ for }\ t\ge 1;\\
\lambda (t) = 0 & \mbox{ for }\ 0 \le t \le 1.
\end{cases}
$$
Then,
$$\lambda (t) \le \lambda_n (t)\le \dfrac 1{\, n\,} \log (n+1) + \lambda (t).$$
Thus, over $[0, \infty), \lambda_n (t)$ converges uniformly to $\lambda (t)$.
For the map $f: \Omega\to {\mathbb C}$,
$$\gamma_n = \dfrac 1{\, n\,} \log (1+ |f|^2 + \cdots + |f|^{2n}) = \lambda_n (|f|^2) ,$$
so that $\gamma_n$ converges uniformly to $\lambda (|f|^2) = \varphi$, and it
follows that $\sqrt{-1} \partial \overline \partial \gamma_n$ converges to
$\sqrt{-1} \partial \overline \partial \varphi$ as positive $(1,1)$-currents on
$\Omega$, as desired. $\quad \square$

\bigskip
\noindent(2.4)
We proceed to give a proof of the Main Theorem.

\bigbreak\noindent {\it Proof.} \   In the language of probability
theory,
$$\left(a_{j_1,\cdots,j_\nu}\right)_{0\leq \nu \leq n, 1\leq j_1\leq\ell, \cdots, 1\leq j_{\nu}\leq
\ell}$$ are independent complex Gaussian random variables of mean
$0$ and variance $1$.  Let $N_{\ell,n}$ be the number of elements in
$$\left(a_{j_1,\cdots,j_\nu}\right)_{0\leq \nu \leq n, 1\leq j_1\leq\ell,\cdots,1\leq j_{\nu}\leq
\ell},$$ which is
$$
N_{\ell, n}=1+\ell+\ell^2+\cdots+\ell^n.
$$
Let $a_0$ be the single indeterminate $a_{j_1,\cdots,j_\nu}$ when
$\nu=0$. By Cauchy's integral formula (or the
Poincar\'e-Lelong formula)
$$
\frac{1}{n}\sum_{G_n(z)=0,\atop z\in\Omega} \delta_z
=\frac{\sqrt{-1}}{n\pi}
\partial\bar\partial\log\left|G_n(z)\right|\leqno{(*)}
$$
on $\Omega$, where $\delta_z$ is the Dirac delta on ${\mathbb C^\ell}$
at the point $z$ of ${\mathbb C^\ell}$. We now consider the normalized
counting divisor ${\mathbf Z}\left(G_n\right)$ of $G_n(z)$ on
$\Omega$ (in the sense of distribution) which is given by
$${\mathbf Z}\left(G_n\right)=
\frac{1}{n}\sum_{G_n(z)=0,\atop z\in\Omega}\delta_z.
$$
By $(*)$, the expectation ${\mathbf E}\left({\mathbf
Z}\left(G_n\right)\right)$ of ${\mathbf Z}\left(G_n\right)$ is equal
to
$$
\displaylines{\qquad\int_{\left(a_{j_1,\cdots,j_\nu}\right)\in{\mathbb
C}^{N_{\ell,n}}} \left(\frac{\sqrt{-1}}{n\pi}
\partial\bar\partial\log\left|G_n(z)\right|\right)\hfill\cr\hfill\times\prod_{\left(a_{j_1,\cdots,j_\nu}\right)\in{\mathbb
C}^{N_{\ell,n}}}\left(\frac{1}{\pi}\,e^{-\left|a_{j_1\cdots
j_\nu}\right|^2}\frac{\sqrt{-1}}{2} da_{j_1\cdots j_\nu}\wedge
d\overline{a_{j_1\cdots j_\nu}}\right).\qquad\cr}
$$
We introduce two column vectors
$$
\vec{\mathbf a}=\left[a_{j_1,\cdots,j_\nu}\right]_{0\leq\nu\leq n,
1\leq j_1\leq\ell,\cdots,1\leq j_\nu\leq \ell}
$$ and
$$\vec {\mathbf
v}(z)=\left[g_\nu(z)f_{j_1}(z)\cdots
f_{j_\nu}(z)\right]_{0\leq\nu\leq n, 1\leq j_1\leq\ell,\cdots,1\leq
j_\nu\leq \ell}
$$
of $N_{\ell,n}$ components each. Here we set $f_0(z)=1$. Then
$G_n(z)$ is equal to the inner product
$$
\left<\vec{\mathbf a},\,\vec{\mathbf v}(z)\right>=\sum_{0\leq\nu\leq
n, 1\leq j_1\leq\ell,\cdots,1\leq j_\nu\leq
\ell}a_{j_1,\cdots,j_\nu}g_\nu(z)f_{j_1}(z)\cdots f_{j_\nu}(z)
$$
of the two $N_{\ell,n}$-vectors $\vec{\mathbf a}$ and $\vec{\mathbf
v}(z)$. The length of the $N_{\ell,n}$-vector $\vec{\mathbf v}(z)$
is given by
$$\left\|\vec{\mathbf
v}(z)\right\|=\left(\sum_{0\leq\nu\leq n, 1\leq
j_1\leq\ell,\cdots,1\leq j_\nu\leq
\ell}|g_\nu(z)|^2|f_{j_1}(z)|^2\cdots|f_{j_\nu}(z)|^2\right)^{\frac{1}{2}}.$$
Introduce the unit $N_{\ell,n}$-vector
$$
\displaylines{\vec{\mathbf u}(z)=\frac{1}{\left\|\vec{\mathbf
v}(z)\right\|}\,{\mathbf v}(z)\cr=\frac{1}{\left(\sum_{0\leq\nu\leq
n, 1\leq j_1, \cdots, j_\nu\leq
\ell}|g_\nu(z)|^2|f_{j_1}(z)|^2\cdots|f_{j_\nu}(z)|^2\right)^{\frac{1}{2}}}\left[g_\nu(z)f_{j_1}(z)\cdots
f_{j_\nu}(z)\right]_{0\leq\nu\leq n, 1\leq j_1, \cdots, j_\nu\leq
\ell}\cr}
$$
in the same direction as $\vec{\mathbf v}(z)$.  Then
$$
\displaylines{\log\left|G_n(z)\right|=\log\left|\left<\vec{\mathbf
a},\,\vec{\mathbf v}(z)\right>\right|=\log\left|\left<\vec{\mathbf
a},\,\left\|\vec{\mathbf v}(z)\right\|\vec{\mathbf
u}(z)\right>\right|\cr =\log\left\|\vec{\mathbf
v}(z)\right\|+\log\left|\left<\vec{\mathbf a},\,\vec{\mathbf
u}(z)\right>\right|.\cr}
$$
Now ${\mathbf E}\left({\mathbf Z}\left(G_n\right)\right)$ is equal
to
$$
\displaylines{\qquad\qquad\int_{\left(a_{j_1,\cdots,j_\nu}\right)\in{\mathbb
C}^{N_{\ell,n}}}
\left(\frac{\sqrt{-1}}{n\pi}\partial\bar\partial\left(\log\left\|\vec{\mathbf
v}(z)\right\|+\log\left|\left<\vec{\mathbf a},\,\vec{\mathbf
u}(z)\right>\right|\right)\right)\hfill\cr\hfill\times\prod_{\left(a_{j_1,\cdots,j_\nu}\right)\in{\mathbb
C}^{N_{\ell,n}}}\left(\frac{1}{\pi}\,e^{-\left|a_{j_1\cdots
j_\nu}\right|^2}\frac{\sqrt{-1}}{2} da_{j_1\cdots j_\nu}\wedge
d\overline{a_{j_1\cdots j_\nu}}\right).\qquad\cr}
$$
Let $ \vec{{\mathbf e}}_0$ be the $N_{\ell,n}$-vector
$$\left(e_{j_1,\cdots,j_\nu}\right)_{0\leq\nu\leq n, 1\leq j_1\leq\ell,\cdots,1\leq j_\nu\leq
\ell}$$ whose only nonzero component is $e_0=1$.  Here comes the key
point of the whole argument. For fixed $z$, we integrate
$$
\displaylines{\qquad\qquad\int_{\left(a_{j_1,\cdots,j_\nu}\right)\in{\mathbb
C}^{N_{\ell,n}}}\left(\frac{\sqrt{-1}}{n\pi}\partial\bar\partial\log\left|\left<\vec{\mathbf
a},\,\vec{\mathbf
u}(z)\right>\right|\right)\hfill\cr\hfill\times\prod_{\left(a_{j_1,\cdots,j_\nu}\right)\in{\mathbb
C}^{N_{\ell,n}}}\left(\frac{1}{\pi}\,e^{-\left|a_{j_1\cdots
j_n}\right|^2}\frac{\sqrt{-1}}{2} da_{j_1\cdots j_n}\wedge
d\overline{a_{j_1\cdots j_n}}\right)\qquad\cr
=\int_{\left(a_{j_1,\cdots,j_\nu}\right)\in{\mathbb
C}^{N_{\ell,n}}}\left(\frac{\sqrt{-1}}{n\pi}\partial\bar\partial\log\left|\left<\vec{\mathbf
a},\,\vec{\mathbf
u}(z)\right>\right|\right)\frac{1}{\pi^{N_{\ell,n}}}\,
e^{-\left\|\vec{\mathbf a}\right\|^2}\cr
=\int_{\left(a_{j_1,\cdots,j_\nu}\right)\in{\mathbb C}^{N_{\ell,n}}}
\left(\frac{\sqrt{-1}}{n\pi}\partial\bar\partial\log\left|\left<\vec{\mathbf
a},\,\vec{{\mathbf e}}_0\right>\right|\right)
\frac{1}{\pi^{N_{\ell,n}}}\, e^{-\left\|\vec{\mathbf
a}\right\|^2}\cr =\int_{\left(a_{j_1,\cdots,j_\nu}\right)\in{\mathbb
C}^{N_{\ell,n}}}\left(\frac{\sqrt{-1}}{n\pi}\partial\bar\partial\log\left|a_0
\right|\right)\frac{1}{\pi^{N_{\ell,n}}}\, e^{-\left\|\vec{\mathbf
a}\right\|^2}\cr=\frac{\sqrt{-1}}{n\pi}\partial\bar\partial\int_{a_0\in{\mathbb
C}}\left(\log\left|a_0\right|\right)\frac{1}{\pi}e^{-\left|a_0\right|^2}\cr}
$$
which is equal to
$$\frac{\sqrt{-1}}{n\pi}\partial\bar\partial A=0\quad{\rm with\ \ }A=\int_{a_0\in{\mathbb
C}}\left(\log\left|a_0\right|\right)\frac{1}{\pi}e^{-\left|a_0\right|^2},$$
because $A$ is a constant.  Note that the equality
$$
\displaylines{\int_{\left(a_{j_1,\cdots,j_\nu}\right)\in{\mathbb
C}^{N_{\ell,n}}}\left(\frac{\sqrt{-1}}{n\pi}\partial\bar\partial\log\left|\left<\vec{\mathbf
a},\,\vec{\mathbf
u}(z)\right>\right|\right)\frac{1}{\pi^{N_{\ell,n}}}\,
e^{-\left\|\vec{\mathbf a}\right\|^2}\cr
=\int_{\left(a_{j_1,\cdots,j_\nu}\right)\in{\mathbb C}^{N_{\ell,n}}}
\left(\frac{\sqrt{-1}}{n\pi}\partial\bar\partial\log\left|\left<\vec{\mathbf
a},\,\vec{{\mathbf e}}_0\right>\right|\right)
\frac{1}{\pi^{N_{\ell,n}}}\, e^{-\left\|\vec{\mathbf
a}\right\|^2}\cr}
$$
in the above string of equalities comes from the fact that for any
fixed $z\in{\mathbb C}$ some unitary transformation of ${\mathbb
C}^{N_{\ell,n}}$ (which may depend on $z$) maps $\vec{\mathbf u}(z)$
to $\vec{{\mathbf e}}_0$ and that $e^{-\left\|\vec{\mathbf
a}\right\|^2}$ is unchanged under any unitary transformation acting
on $\vec{\mathbf a}$.  Thus the limit of ${\mathbf E}\left({\mathbf
Z}\left(G_n\right)\right)$ as $n\to\infty$ is equal to
$$
\lim_{n\to\infty}\int_{\left(a_{j_1,\cdots,j_\nu}\right)\in{\mathbb
C}^{N_{\ell,n}}}\left(\frac{\sqrt{-1}}{n\pi}\partial\bar\partial\log\left\|\vec{\mathbf
v}(z)\right\|\right)\frac{1}{\pi^{N_{\ell,n}}}\,
e^{-\left\|\vec{\mathbf a}\right\|^2},$$ which after integration
over $$\left(a_{j_1,\cdots,j_\nu}\right)_{0\leq\nu\leq n, 1\leq
j_1\leq\ell,\cdots,1\leq j_\nu\leq \ell}$$ is simply equal to
$$
\displaylines{\lim_{n\to\infty}\frac{\sqrt{-1}}{n\pi}
\partial\bar\partial\log\left\|
\vec{\mathbf v}(z)\right\|\cr
=\lim_{n\to\infty}\frac{1}{n}\left(\frac{\sqrt{-1}}{2\pi}\partial\bar\partial\log\sum_{k=0}^n|g_k(z)|^2\left(\sum_{j=1}^\ell
\left|f_j(z)\right|^2\right)^k\right)\cr}
$$
From (2.1) Proposition, we have
$$
\lim_{n\to\infty}\frac{1}{n}\left(\frac{\sqrt{-1}}{2\pi}\partial\bar\partial\log\sum_{k=0}^n\left(\sum_{j=1}^\ell
\left|f_j(z)\right|^2\right)^k\right)
$$
is equal to  the pullback by $f$ of
$$\left|w\right|\Xi\left(\left|w\right|\right)\frac{\sqrt{-1}}{2\pi}\partial\bar\partial\log\left|w\right|+
\left[\delta_{{}_{S^{2\ell-1}}}\right]\wedge
\frac{1}{\left|w\right|^2}\frac{\sqrt{-1}}{2}\sum_{j=1}^\ell\left(w_j
d\bar w_j-\bar w_j dw_j\right),
$$
where $w\in{\mathbb C}^\ell=\left(w_1,\cdots,w_\ell\right)$ is
variable in the target space of the map
$f=\left(f_1,\cdots,f_\ell\right):{\mathbb C}\to{\mathbb C}^\ell$. By computation
$$
\displaylines{\lim_{n\to\infty}\frac{1}{n}\frac{\sqrt{-1}}{2\pi}\partial\bar\partial\log\left(\sum_{k=0}^n|g_k(z)|^2\left(\sum_{j=1}^\ell
\left|f_j(z)\right|^2\right)^k\right)\cr
-\lim_{n\to\infty}\frac{1}{n}\frac{\sqrt{-1}}{2\pi}\partial\bar\partial\log\left(\sum_{k=0}^n\left(\sum_{j=1}^\ell
\left|f_j(z)\right|^2\right)^k\right)\cr }
$$
$$
=\lim_{n\to\infty}\frac{1}{n}\frac{\sqrt{-1}}{2\pi}\partial\bar\partial\left(\log\sum_{k=0}^n|g_k(z)|^2\left(\sum_{j=1}^\ell
\left|f_j(z)\right|^2\right)^k-\log\sum_{k=0}^n\left(\sum_{j=1}^\ell
\left|f_j(z)\right|^2\right)^k\right)
$$
$$
=\frac{\sqrt{-1}}{2\pi}\partial\bar\partial\lim_{n\to\infty}\frac{1}{n}\left(\log\sum_{k=0}^n|g_k(z)|^2\left(\sum_{j=1}^\ell
\left|f_j(z)\right|^2\right)^k-\log\sum_{k=0}^n\left(\sum_{j=1}^\ell
\left|f_j(z)\right|^2\right)^k\right)
$$
$$
=\frac{\sqrt{-1}}{2\pi}\partial\bar\partial\lim_{n\to\infty}\frac{1}{n}\log\left(\frac{\sum_{k=0}^n|g_k(z)|^2\left(\sum_{j=1}^\ell
\left|f_j(z)\right|^2\right)^k}{\sum_{k=0}^n\left(\sum_{j=1}^\ell
\left|f_j(z)\right|^2\right)^k}\right)\, .
$$
Without loss of generality, we may assume that for any $z\in\Omega$, $A_z$, $B_z\ge 1$.
Granted this, replacing $\kappa_{z,n}$ by $\max\big(\kappa_{z,0},\ldots,\kappa_{z,n}\big)$, etc., without
loss of generality we may assume that the four sequences $\kappa_{z,n}$,
$\lambda_{z,n}$, $\xi_{z,n}$ and $\eta_{z,n}$ are non-decreasing sequences.
By (iii) for every $z\in \Omega$ there exists a positive constant $c_z$ and
a positive integer $J(z)$ such that whenever $j\ge J (z)$ we have
$$
A_z^{\xi_{z,j}} \big(1+|f(z)|\big)^{\eta_{z,j}} |g_j(z)| \ge c_z.
$$
(Here and in what follows to streamline the notations we will write
$A_z^{\xi_{z,j}}$ to mean $(A_z)^{\xi_{z,j}}$, etc.)
For every $z\in \Omega$ we have
$$\begin{array}{ll}
& \sum_{k=0}^n |g_k(z)|^2 |f(z)|^{2k} \ge \max\big(|g_0(z)|^2, |g_n(z)|^2 |f(z)|^{2n}\big) \\
\ge & \max \big(|g_0(z)|^2, c_z^2 A_z^{-2\xi_{z,n}} \big(1+|f(z)|\big)^{-2\eta_{z,n}}
|f(z)|^{2n}\big)\ .
\end{array}
$$
On the other hand, when $|f(z)|\le 1$ we have
$$
\sum_{k=0}^n |g_k(z)|^2 |f(z)|^{2k} \le (n+1) B_z^{2\kappa_{z,n}} \cdot 4^{\lambda_{z,n}}\ ;
\leqno{(\dagger)}
$$
and, when $|f(z)| \ge 1$ we have
$$
\sum_{k=0}^n |g_k(z)|^2 |f(z)|^{2k} \le (n+1)B_z^{2\kappa_{z,n}}
\big(1+|f(z)|\big)^{2\lambda_{z,n}} |f(z)|^{2n}\ ,
$$
so that
$$
\begin{array}{ll}
\sum_{k=0}^n |g_k(z)|^2 |f(z)|^{2k}
& \le \max \Big((n+1) B_z^{2\kappa_{z,n}} \cdot 4^{\lambda_{z,n}}\ ,\\
& \qquad (n+1)B_z^{2\kappa_{z,n}} \big(1+|f(z)|\big)^{2\lambda_{z,n}} |f(z)|^{2n}\Big)\ .
\end{array}
$$
Similarly for the function $\Big(\sum\limits_{k=0}^n |f(z)|^{2k}\Big)^{\frac 1n}$ we have
$$
\max\big(1, |f(z)|^2\big) \le \Big(\sum_{k=0}^n |f(z)|^{2k}\Big)^{\frac 1n}
\le (n+1)^{\frac 1n} \max \big(1, |f(z)|^2\big)\ .
$$
Finally, recalling that
$$
\big(h_n(z)\big)^{\frac 1n} =
\left(\frac{\sum\limits_{k=0}^n |g_k(z)|^2 |f(z)|^{2k}}{\sum\limits_{k=0}^n
|f(z)|^{2k}}\right)^{\frac 1n}\ ,
$$
we have, for $z\in \Omega$,
$$
\frac{\max\big(|g_0(z)|^2, c_z^2 A_z^{-2\xi_{z,n}} \big(1+|f(z)|\big)^{-2\eta_{z,n}}\big)^{\frac 1n}} {(n+1)^{\frac 1n} \max\big(1,|f(z)|^2\big)}
\le \big(h_n(z)\big)^{\frac 1n}
$$
$$
\le \frac{\max\big((n+1)B_z^{2\kappa_{z,n}}\cdot 4^{\lambda_{z,n}},(n+1)B_z^{2\kappa_{z,n}} \big(1+|f(z)|\big)^{2\lambda_{z,n}} |f(z)|^{2n}\big)^{\frac 1n}} {\max\big(1,|f(z)|^2\big)}\ .
$$
For the lower bound of $\big(h(z)\big)^{\frac 1n}$ we note that
$$
\lim\limits_{n\to\infty} |g_0(z)|^{\frac 2n} = 1\ ;
$$
$$
\lim\limits_{n\to\infty} \big(c_z^2A_z^{-2\xi_{z,n}} \big(1+|f(z)|\big)^{-2\eta_{z,n}}\big)^{\frac 1n} \big|f(z)\big|^{2n} $$
$$
= \lim\limits_{n\to\infty} c_z^{\frac 2{\,n\,}} A_z^{-\frac{2\xi_{z,n}} n} \big(1+|f(z)|\big)^{-\frac{2\eta_{z,n}}n} \big|f(z)\big|^2 = \big|f(z)\big|^2
$$
where we have used the assumptions
$\lim\limits_{n\to\infty} \frac{\xi_{z,n}}n = \lim\limits_{n\to\infty} \frac{\eta_{z,n}}n=0$.
For the upper bound of $\big(h_n(z)\big)^{\frac 1n}$
we note that
$$
\lim\limits_{n\to\infty} \big((n+1)B_z^{2\kappa_{z,n}} \cdot 4^{\lambda_{z,n}}\big)^{\frac 1n}
= \lim\limits_{n\to\infty} (n+1)^{\frac 1n} B_z^{\frac{2\kappa_{z,n}}n} 4^{\frac{\lambda_{z,n}}n} = 1\ ;
$$
$$
\lim\limits_{n\to\infty} \big((n+1)B_z^{2\kappa_{z,n}} \big(1+|f(z)|\big)^{2\lambda_{z,n}}
\big|f(z)\big|^{2n}\big)^{\frac 1n}
$$
$$
= \lim\limits_{n\to\infty} (n+1)^{\frac 1n} B_z^{\frac{2\kappa_{z,n}}n}
\big(1+|f(z)|\big)^{\frac{2\lambda_{z,n}}n} \big|f(z)\big|^2 = \big|f(z)\big|^2
$$
where we have used the assumptions
$\lim\limits_{n\to\infty} \frac{\kappa_{z,n}}n = \lim\limits_{n\to\infty}
\frac{\lambda_{z,n}}n = 0$.
Thus, for any $z\in \Omega$ we have
$$
1 = \frac{\max\big(1,|f(z)|^2\big)}{\max\big(1,|f(z)|^2\big)}
\le \varliminf\limits_{n\to\infty}h_n(z)^{\frac 1{\,n\,}}
\le \varlimsup\limits_{n\to\infty}h_n(z)
\le \frac{\max\big(1,|f(z)|^2\big)}{\max\big(1,|f(z)|^2\big)} = 1\ .
$$
so that
$$
\lim\limits_{n\to\infty} h_n(z)^{\frac 1n} = 1\ ;\quad
\lim\limits_{n\to\infty} \log \big(h_n(z)^{\frac 1{\,n\,}}\big) = 0\ .
$$

\smallskip\noindent
Under the assumptions of the Main Theorem write
$$
\varphi_n = \log \Big(\sum_{j=0}^n \big|g_j(z)\big|^2 \big|f(z)\big|^{2j}\Big)^{\frac 1n}\ .
$$
Then, $\log h_n{}^{\frac 1{\,n\,}} = \varphi_n - \gamma_n$. Since
$\gamma_n$ converges to $\varphi = \log^+ |f|^2$ by Lemma 1 and
$\log h_n{}^{\frac 1{\,n\,}} $ converges pointwise to 0,
we conclude that $\varphi_n(z)$ converges to $\varphi (z)$
for every $z\in \Omega$. Clearly
$\varphi_n$ and $\varphi$ are continuous subharmonic
functions on $\Omega$. Moreover from
$(\dagger)$ we have for every $z\in \Omega$
$$
\varphi_n (z) \le \dfrac 1{\,n\,} \log (n+1)
+ \dfrac{2\kappa_{z,n}} n \log B (z) +
\dfrac{\lambda_{z,n}} n \log 4,
$$
and by assumption on any compact subset $K \subset \Omega$, $B_z$ and
the sequence of functions $\dfrac{\kappa_{z,n}} n$ and
$\dfrac{\lambda_{z,n}} n$ are uniformly bounded from above by some constant
$c_K$ for $z\in K$, and we conclude that the sequence of subharmonic
functions $\big ( \varphi_n(z) \big)^\infty_{n=0}$ are
uniformly bounded from above on compact subsets. Finally, we make use of
Lemma 2 below on the convergence of positive $(1,1)$ currents.
Granting Lemma 2, the Main Theorem follows readily.
\quad $\square$

\smallskip
The discussion below involves distributions on a domain in $\mathbb C$.
Denote by $d\lambda$ the Lebesgue measure on $\mathbb C$.
Any locally integrable function $s$ on $\Omega$ defines a distribution $T_s$
on $\Omega$ given by $T_s(\rho) = \displaystyle \int_{\Omega} s\rho \ d\lambda$
for any smooth function $\rho$ on $\Omega$ of compact support,
and in what follows we will identify $s$ with the distribution $T_s$ it defines.
There is a standard process for smoothing distributions, as follows.
Let $\chi$ be a nonnegative smooth function on $\mathbb C$ of support
lying on the unit disk $\Delta$ such that $\chi(e^{i\theta}z) = \chi(z)$ for any $z \in \mathbb C$
and any $\theta \in \mathbb R$, and for any $\epsilon > 0$ write $\chi_{\epsilon}(z) =
\chi\Big (\frac {\phantom{,}z\phantom{,}}{\epsilon}\Big )$.
For a distribution $Q$ defined on some domain in $\mathbb C$ and for $\epsilon > 0$
we write $Q_{\epsilon} := Q * \chi_{\epsilon}$
wherever the convolution is defined.  We have the following elementary
lemma on positive currents associated to subharmonic functions, for which a proof
is included below for easy reference.

\smallskip\noindent
{\bf Lemma 2}\
{\it
Let $\Omega\subset\Bbb C$ be a plane domain.
Suppose $(\varphi_n)_{n=0}^\infty$ is a sequence of subharmonic functions on $\Omega$ such that
$\varphi_n(z)$ are uniformly bounded from above on each compact subset $K$ of $\Omega$.
Assume that $\varphi_n$ converges pointwise to some continuous {\rm(}subharmonic{\rm)} function $\varphi$.
Then,
$\lim\limits_{n\to\infty}\varphi_n = \varphi$ in $L^1_{\text{loc}}(\Omega)$.  As a consequence,
$\sqrt{-1}\,\partial\bar\partial\varphi_n$ converges to $\sqrt{-1}\,\partial\bar\partial\varphi$
in the sense of currents.
}

\smallskip\noindent
{\it Proof of Lemma 2.}
Let $D = \Delta (a; r)$ be any disk centred at $a\in \Omega$ of radius
$r >0$ such that $\overline D \subset \Omega$.
We claim that the Lebesgue integrals
$\displaystyle \int_{\Delta(a;r)} |\varphi_n| \ d\lambda$ are bounded independent of $n$.
Without loss of generality, we may assume that $\varphi \le 0$ on $\overline D$.  By the
sub-mean-value inequality for subharmonic functions we have
$$
\varphi_n(a) \le \frac 1{\pi r^2}\int_{\Delta(a;r)} \varphi_n(\zeta) \ d\xi \ d\eta
$$
where $\zeta = \xi + \sqrt{-1}\eta$ is the Euclidean coordinate of the variable
of integration $\zeta$, showing that the integral of $-\varphi_n$
over $\Delta(a;r)$ are bounded independent of $n$.  Covering $\Omega$ by a
countable and locally finite family of relatively compact open disks $D$,
it follows that on any compact subset $K \subset \Omega$ the $L^1$-norms of
$\varphi_n|_K$ are bounded independent of $n$.  As a consequence,
given any subsequence
$\varphi_{\sigma(n)}$ of $\varphi_n$, some subsequence $\psi_n := \varphi_{\sigma(\tau(n))}$
of $\varphi_{\sigma(n)}$ must converge to a distribution
$S$ on $\Omega$. We claim that any such a limit must be given by the (continuous)
subharmonic function $\varphi$.
As a consequence, $\varphi_n$ converges to $\varphi$ in
$L^1_{\rm loc}(\Omega)$.

Since $\psi_n$ converges to the distribution $S$, for any
$\epsilon > 0$, $\varphi_{n,\epsilon}$ converges to the smooth function $S_{\epsilon}$ as
$n$ tends to $\infty$.  Since $\psi_n$ is subharmonic, $\psi_{n,\epsilon}$ is monotonically decreasing
as $\epsilon \mapsto 0$ for each nonnegative integer $n$, and it follows readily that
$S_{\epsilon}$ is also monotonically decreasing as $\epsilon \mapsto 0$.
Hence, $S$ is the limit as a distribution of the smooth functions $S_{\epsilon}$. Writing $\psi(z)
:= \lim_{\epsilon \mapsto 0}S_{\epsilon}(z)$, by the Monotone Convergence Theorem the distribution
$S$ is nothing other than the function $\psi$,
which is in particular locally integrable.  Since $\psi_n$
converges to $S$ as distributions, we conclude that $\varphi_{\sigma(n)} = \psi_n$ converges to $\psi$ in
$L^1_{\rm loc}(\Omega)$, implying that $\psi_n$ converges pointwise to $\psi$
almost everywhere on $D$.  However, by assumption $\psi_n = \varphi_{\sigma(n)}$
converges pointwise to $\varphi$, hence $\varphi$ and $\psi$ must agree almost everywhere
on $\Omega$. In particular, $\varphi_n$ must converge to $\varphi$ in $L^1_{\rm loc}(\Omega)$.
The proof of Lemma 2 is complete.  $\quad \square$

\smallskip\noindent
\section{Acknowledgments}

\smallskip\noindent
The author would like to thank Professors Yum-Tong Siu and
Ngaiming Mok for discussions and suggestions concerning this article.

\section{References}
\medbreak\noindent[DSZ-2004] Michael Douglas, Bernard Shiffman, and
Steve Zelditch, Critical points and supersymmetric vacua. I. {\it
Comm. Math. Phys.} 252 (2004), no. 1-3, 325--358.

\medbreak\noindent[DSZ-2006a] Michael R. Douglas, Bernard Shiffman,
and Steve Zelditch, Critical points and supersymmetric vacua. II.
Asymptotic and extremal metrics. {\it J. Differential Geom.} 72
(2006), no. 3, 381--427.

\medbreak\noindent[DSZ-2006b] Michael Douglas, Bernard Shiffman, and
Steve Zelditch, Critical points and supersymmetric vacua. III.
String/M models. {\it Comm. Math. Phys.} 265 (2006), no. 3,
617--671.

\medbreak\noindent[H] J M. Hammersley,  The zeros of a random
polynomial. In: Proceedings of the Third Berkeley Symposium on
Mathematical Statistics and Probability, 1954¨C1955, vol. II,
89¨C111, Berkeley-Los Angeles: University of California Press, 1956

\medbreak\noindent[SV] L.A. Shepp and R.J. Vanderbei, The complex
zeros of random polynomials. {\it Trans. Am. Math. Soc.}
\textbf{347} (1995), 4365--4384.

\medbreak\noindent[SZ-1999] Bernard Shiffman and Steve Zelditch,
Distribution of zeros of random and quantum chaotic sections of
positive line bundles. {\it Comm. Math. Phys.} 200 (1999), no. 3,
661--683.

\medbreak\noindent[Y-2008] W. Yao, Distribution of
Normalized Zero-Sets of Random Entire Functions,
arxiv.org/abs/0811.3365.

\end{document}